\documentclass[nolineno]{biometrika}

\usepackage{amsmath}

\usepackage{times}
\usepackage{bm}
\usepackage{natbib}
\usepackage{graphicx}
\usepackage{color}

\usepackage[plain,noend]{algorithm2e}
\def\T{{ \mathrm{\scriptscriptstyle T} }}

\def\ft{\mathrm{Ft}}
\def\Rebu{\mathrm{Re}}

\newcommand{\var}{\mbox{var}}

\def\JRSSB{{\sl Journal of the Royal Statistical Society}, {\bf B}}

\def\BKA{{\sl Biometrika}}
\def\JASA{{\sl Journal of the American Statistical Association}}

\def\SS{{\sl Statistica Sinica}}

\def\AS{{\sl The Annals of Statistics}}

\def\JOE{{\sl Journal of Econometrics}}

\def\SPL{{\sl Statistics and Probability Letters}}

\begin{document}

\jname{Biometrika}
%% The year, volume, and number are determined on publication
\jyear{2017}
\jvol{103}
\jnum{1}
%% The \doi{...} and \accessdate commands are used by the production team
%\doi{10.1093/biomet/asm023}
\accessdate{Advance Access publication on 31 July 2016}

%% These dates are usually set by the production team
\received{2 January 2017}
\revised{1 April 2017}

%% The left and right page headers are defined here:
\markboth{J. Chang, A. Delaigle, P. Hall \and C. Y. Tang}{A frequency domain analysis of high frequency data}

%% Here are the title, author names and addresses
\title{A frequency domain analysis of the error distribution from noisy high-frequency data}

\author{JINYUAN CHANG }
\affil{School of Statistics, Southwestern University of Finance and Economics, Chengdu,
Sichuan 611130, China \email{changjinyuan@swufe.edu.cn} }

\author{AURORE DELAIGLE, PETER HALL}
\affil{School of Mathematics and Statistics,  University of Melbourne, Parkville, Victoria 3010, Australia \email{aurored@unimelb.edu.au}}

\author{\and CHENG YONG TANG}
\affil{Department of Statistical Science, Temple University, Philadelphia, Pennsylvania 19122-6083, U.S.A. \email{yongtang@temple.edu}}

\maketitle

\begin{abstract}
Data observed at high sampling frequency are typically assumed to be an additive composite of a relatively slow-varying continuous-time component,  a latent stochastic process or a  smooth random function,  and measurement error.
Supposing that the latent component is an It\^{o} diffusion process,
we propose to estimate the measurement error density function by applying a deconvolution technique with appropriate localization.
Our estimator, which does not require equally-spaced observed times, is consistent and minimax rate optimal.  We also investigate estimators of the moments of the error distribution and their properties,
propose a frequency domain estimator for the integrated volatility of the underlying stochastic process, and show that it achieves the optimal convergence rate.
Simulations and a real data analysis validate our analysis.
\end{abstract}


\section{Introduction}\label{sec:intro}

High-frequency data, observed sequentially at small time intervals, arise in various settings and applications. For example, in  social and behavioural investigations, such data are often collected  in so-called  intensive longitudinal studies; see \cite{BL13}. In functional data analysis,  the observations are often considered to be of one of two types: the dense setting, which corresponds to high-frequency data, and the sparse setting, where the data are of low frequency; see, among others,   \cite{Wang_2016} and  \cite{MSS2011}.
In finance, analyzing high-frequency intra-day transaction data has received increasing attention; see  \cite{Hautsch_2012}, \cite{AJ_2014_Book}  and the references therein.

High-frequency  observations are often contaminated by measurement errors. For example, in the dense functional data setting, it is common to assume that the observed  discrete data are a noisy version  of an underlying unknown smooth curve. In finance too, high-frequency data are often regarded as a noisy version of a latent continuous-time stochastic process, observed at consecutive discrete time points. The latent process is considered to have a continuous path, and the measurement error represents the market microstructure noise; see  \cite{AitYu_2009_AOAS}.
Since increasing the sampling frequency implies smaller sampling errors caused by the discretization  of  the underlying continuous-time process, we could expect that high-frequency data enable more accurate inference.  However, as the sampling frequency increases, the difference between nearby observations is dominated by random noise, which makes standard methods of inference inconsistent; see, for example,  \cite{Zhangetal_2005_JASA}. Therefore,  a main concern in high-frequency financial data analysis has been to find ways of recovering the signal of quantities of interest from noisy high-frequency observations;  see \cite{AJ_2014_Book}.

Revealing the distributional properties of the measurement errors is crucial for recovering the signal of the continuous-time  process from noisy high-frequency  data.
First,  many estimation procedures
require distributional assumptions  on the measurement errors. Second,
statistical inference including hypothesis testing and confidence set estimation  inevitably involves unknown nuisance parameters  determined by the distribution of the measurement errors;  see, for example, \cite{Zhang_2006_Bern}, \cite{Xiu_2010_JOE},   \cite{Xiuetal_2010_JASA},  and \cite{LiuTang_2014_JOE}.

The measurement errors themselves sometimes contain  useful information, both  theoretical and practical, so that successfully recovering the measurement error distribution can improve our understanding of  data structures. For example, \cite{Duf_2014} argue that the microstructure noise can help understand financial crises, and \cite{AitYu_2009_AOAS} make a connection between the microstructure noise  and the liquidity of stocks. In longitudinal studies, the measurement errors can help reveal interesting characteristics of a population, such as the way in which individuals misreport their diet.
\cite{Jacodetal_2015} recently  highlighted the importance of  statistical properties of the measurement errors, and  studied the estimation of moments.
Despite these important developments, to our knowledge, no method has been proposed for estimating the entire distribution of the errors.

In this paper, we consider frequency domain analyses of high-frequency  data, focusing on the measurement errors contaminating the continuous-time processes.
In high-frequency observations,  the relative magnitude of the changes in values of the underlying continuous-time  process is small. Compared with the measurement errors, it becomes negligible,  locally in a small neighborhood of a given time point.   As a result, estimating the error distribution shares some  common features with nonparametric measurement error problems with repeated measurements studied in  \cite{DelaigleHallMeister_2008}, and with nonparametric estimation from aggregated data studied by \cite{Meister_2007} and \cite{DelaigleZhou}, where the estimation techniques require working in the Fourier domain.

Motivated by this,   we propose to estimate the characteristic function of the measurement errors by locally averaging the empirical characteristic functions of the changes in values of the high-frequency data. We obtain a nonparametric estimator of the probability density function of the measurement errors and show  that it is consistent and minimax rate optimal.
We propose a simple method for consistently estimating the moments of the measurement errors.
Using our estimator of the characteristic function of the errors,  we develop a new rate-optimal multiscale frequency domain estimator of the integrated volatility of the stochastic process, a key quantity of interest in high-frequency financial data analysis.

In the sequel, for two sequences of real numbers $\{a_n\}$ and $\{b_n\}$, we write $a_n\asymp b_n$ if there exist positive constants $c_1$ and $c_2$ such that $c_1\leq a_n/b_n\leq c_2$ for all $n\geq1$. We denote by $f*g$ the convolution of two functions $f$ and $g$ defined by $f*g(s)=\int f(s-\tau)g(\tau)\,d\tau$.

\section{Methodology}
\label{s2}

\subsection{Model and data}\label{sec:model}
We are interested in a continuous-time process $(X_t)_{t\in[0,T]}$ observed at high frequency with  $T>0$.
We assume that $X_t$ follows a diffusion process
\begin{equation}\label{eq:model}
dX_t=\mu_t~dt+\sigma_t~dB_t\,,
\end{equation}
where the drift $\mu_t$ is a locally bounded and progressively measurable process, $\sigma_t$ is a positive and locally bounded It\^{o} semimartingale, and $B_t$ is a standard Brownian motion. The process $\sigma_t^2$ represents the volatility of the process $X_t$ at time $t$, and is often investigated in its integrated form, $\int_0^T \sigma_t^2\,dt$, called integrated volatility.

\begin{remark}
The model at \eqref{eq:model} is often used when $X_t=\log S_t$, where $S_t$ denotes the price process of an equity; see for example \cite{Zhangetal_2005_JASA}. It is also used to model applications in biology, physics and many other fields; see \cite{OSP_2019}. All theoretical properties of our work are derived under this model, but the methods derived in Sections~\ref{sec:methfU}  and \ref{se:mome} could be applied to other types of processes $X_t$, such as the smooth ones typically encountered in the functional data literature. The key property that makes our methods consistent is the continuity of the underlying process $X_t$, but the convergence rates of our estimators depend on more specific assumptions, such as those implied by the model at \eqref{eq:model}.
\end{remark}

 Our data are observed at a generic discrete grid  $\mathcal {G}=\{t_0,\ldots,t_n\}$ of time points where, without loss of generality, we let $t_0=0$ and $t_n=T$.
The observed data are contaminated by additive measurement errors, so that what we observe is a sample $\{Y_{t_j}\}_{j=0}^n$, where
\begin{equation}\label{eq:finalmodel}
Y_{t_j}=X_{t_j}+U_{t_j}.
\end{equation}
A conventional assumption when analyzing noisy high-frequency  data is that the random measurement error $U_t$ is independent of $X_t$; see \cite{Xiu_2010_JOE},   \cite{Xiuetal_2010_JASA},  and \cite{LiuTang_2014_JOE}. This is also a standard assumption in the measurement error and the functional data literature; see \cite{CH1988}, \cite{SC1990} and \cite{Wang_2016}. Likewise, we make the following assumption:
 \begin{assumption}\label{as:inde}
The errors $\{U_{t_j}\}_{j=0}^n$ are independently and identically distributed, with unknown density $f_U$, and are independent of the process $(X_t)_{t\in[0,T]}$.
\end{assumption}

We are interested in deriving statistical properties of the noise term $U_{t_j}$ when $T$ is fixed and the frequency of observations increases, that is, when $\max_{1\leq j\leq n}\Delta t_j\rightarrow0$ as
$n\rightarrow\infty$, where $\Delta t_j=t_j-t_{j-1}$ for
$j=1,\ldots,n$. Here, the time points $t_j$'s do not need to be equispaced. Formally, we make the following assumption:
\begin{assumption}\label{as:space}
As $n\rightarrow\infty$, $\min_{1\leq j\leq n}\Delta t_j/\max_{1\leq j\leq n}\Delta t_j$ is uniformly bounded away from zero and infinity.
\end{assumption}

Throughout we use $f^{\ft}$ to denote the Fourier transform of a function $f$ and we make the following assumption on the characteristic function of the errors, which is almost always assumed in the related nonparametric measurement error literature:
\begin{assumption}\label{as:symm}
$f_U^{\ft}$ is real-valued and does not vanish at any point on the
real line.
\end{assumption}

\subsection{Estimating the error density $f_U$}\label{sec:methfU}
Motivated by our discussion in Section~\ref{sec:intro}, we wish to estimate the error density $f_U$. At a given $t_j$, if  we had access to repeated noisy measurements of $X_{t_j}$, say $Y_{t_j,k}=X_{t_j}+U_{t_j,k}$, for $k=1,\ldots,r$, where the $U_{t_j,k}$'s are independent and each $U_{t_j,k}\sim f_U$, then for $\ell\neq k$ we would have $(Y_{t_j,\ell}-Y_{t_j,k})\sim f_U*f_U$. Under Assumption~\ref{as:symm}, using the approach in \cite{DelaigleHallMeister_2008}  we could estimate $f_U^{\ft}$ by the square root of the empirical characteristic function of the  $(Y_{t_j,\ell}-Y_{t_j,k})$'s; then by Fourier inversion we could deduce an estimator of~$f_U$.

However, for high-frequency  data, at each given  $t_j$ we have access to only one contaminated measurement $Y_{t_j}$. Therefore, the above  technique cannot be applied.  But since $(X_t)_{t\in[0,T]}$ is a continuous-time and continuous-path
stochastic process, $|X_{t+h}-X_t|\rightarrow0$ almost surely as $h\rightarrow0$.
Thus, the   collection of observations $\{Y_{t_\ell}\}$, where $t_\ell$ lies in a small neighborhood $\cal N$ of $t_j$,  can be approximately viewed as repeated measurements  of $X_{t_j}$ contaminated by independently and identically distributed errors $\{U_{t_\ell}\}$. As the sampling frequency increases,   we have multiple observations in increasingly smaller neighborhoods $\cal N$, which suggests that the density of  $Y_{t_\ell}-Y_{t_j}$, for $t_j\neq t_\ell \in \cal N$,  gets increasingly closer to $f_U*f_U$. Therefore, we can expect that as the sample frequency increases, the approach suggested by \cite{DelaigleHallMeister_2008}, applied to the $Y_{t_\ell}-Y_{t_j}$'s, for $t_\ell$ and $t_j$ sufficiently close, can provide an increasingly accurate estimator of $f_U$.

We shall prove in Section~\ref{sec:theoryfU} that this heuristic is correct as long as the $t_\ell$'s and $t_j$'s are carefully chosen, which we characterize through a distance $\xi$.
For $\mathcal {G}$ defined in Section~\ref{sec:model} and $\xi>0$,  we define
\begin{equation}\label{eq:sj}
S_j=\{t_\ell\in \mathcal {G}:|t_\ell-t_j|\leq \xi~\textrm{and}~\ell\neq j\}~~~~~~(j=0,\ldots,n) \,,
\end{equation}
and denote by $N_j$ the number of points in $S_j$. For a fixed $T$,  Assumption~\ref{as:space} implies that $\min_{1\leq j\leq n}\Delta t_j\asymp\max_{1\leq j\leq n}\Delta t_j\asymp n^{-1}$, so that $\max_{1\leq j\leq n}N_j\asymp\min_{1\leq j\leq n}N_j\asymp n\xi$.

Following the discussion above and recalling Assumption~3, for a given $\xi$ we define our estimator of $f_U^{\ft}(s)$ by the square root of the real part of the empirical characteristic function of the difference of nearby $Y_{t_{j}}$'s:
\begin{equation}\label{eq:char0}
\hat{f}_{U,1}^{\ft}(s;\xi)=\bigg|\frac{1}{N(\xi)}\sum_{j=0}^n\sum_{t_\ell\in S_j}\cos\{s(Y_{t_{\ell}}-Y_{t_{j}})\}\bigg|^{1/2}\,,
\end{equation}
where $N(\xi)=\sum_{j=0}^n N_j$.
Here $\xi$ can be viewed as a parameter controlling the trade-off between the bias and the variance of $\hat{f}_{U,1}^{\ft}$: a smaller $\xi$ results in a smaller bias, but also results  in a smaller $N(\xi)$ so that the variance is higher. On the other hand, a larger $\xi$ induces a lower variance, but  comes at the price of a larger bias due to the contribution from the dynamics in $X_t$. The choice of $\xi$ in practice will be discussed in Section~\ref{ParamSelec}.

It follows from the Fourier inversion theorem that $f_U(x)=(2\pi)^{-1}\int e^{-{i}sx}f_{U}^{\ft}(s)\,ds$, where $i^2=-1$. We can obtain an estimator of $f_U$ by replacing $f_{U}^{\ft}(s)$ with  $\hat{f}_{U,1}^{\ft}(s,\xi)$ in this integral. However, since $\hat{f}_{U,1}^{\ft}(s,\xi)$ is an empirical characteristic function, it is unreliable when $|s|$ is large. In order for the integral to exist, $\hat{f}_{U,1}^{\ft}(s,\xi)$ needs to be multiplied by a regularizing factor that puts less weight on large $|s|$. As the sample size increases, $\hat{f}_{U,1}^{\ft}(s,\xi)$ becomes more reliable, and this should be reflected by letting the weight depend on the sample size. Using standard kernel smoothing techniques, this can be implemented by taking
\begin{equation*}\label{eq:char1}
\hat{f}_{U,2}^{\ft}(s;\xi)=\hat{f}_{U,1}^{\ft}(s;\xi)\,\mathcal{K}^{\ft}(sh)\,,
\end{equation*}
where $\mathcal{K}^\ft$ is the Fourier transform of a kernel function $\mathcal{K}$, and $h>0$ is a bandwidth parameter that satisfies $h\rightarrow0$ as $n\rightarrow\infty$. Then, we define our estimator of $f_U(x)$ by
\begin{equation}\label{eq:density}
\hat{f}_U(x;\xi)=\frac{1}{2\pi}\int e^{-{ i}sx}\hat{f}_{U,1}^{\ft}(s;\xi)\mathcal{K}^{\ft}(sh)\,ds=\frac{1}{2\pi}\int e^{-{i}sx}\hat{f}_{U,2}^{\ft}(s;\xi)\,ds\,.
\end{equation}

For the consistency of a kernel density estimator, the kernel function $\mathcal{K}$ needs to be symmetric and integrate to unity. In practice, the estimator is often not very sensitive to the choice of the kernel compared to that of the bandwidth. Popular choices are the Gaussian kernel, i.e. the standard normal density, and the sinc kernel, whose Fourier transform is  $\mathcal{K}^{\ft}(s)=I(|s|\leq 1)$, where $I(\cdot)$ denotes the indicator function. In practice, an advantage of the Gaussian kernel is that it can produce visually attractive smoother estimators than the sinc kernel. The sinc kernel is more advantageous  analytically; see Section~\ref{sec:theoryfU}.

\subsection{Properties of the density estimator }\label{sec:theoryfU}
Assume that the continuous-time process $(X_t)_{t\in[0,T]}$ in (\ref{eq:model}) belongs to the class $\mathcal{X}(C_1)$ for some $C_1>0$:   %where, for any positive $C_1$, we let
\begin{equation*}\label{eq:xclass}
\begin{split}
\mathcal{X}(C_1)=\bigg\{(X_{t})_{t\in[0,T]}:&~X_t~\textrm{satisfies (\ref{eq:model}) with}\sup_{0\leq t\leq T}E(\mu_t^4)\leq C_1~\textrm{and}\sup_{0\leq t\leq T}E(\sigma_t^4)\leq C_1\bigg\}\,,
\end{split}
\end{equation*}
and that $f_U$ belongs to the class:
\begin{equation*}\label{eq:uclass1}
\begin{split}
\mathcal{F}_1(\alpha,C_2)=&~\big\{f:~ |f^{\ft}(s)|\leq C_2(1+|s|)^{-\alpha}~\textrm{for all real}~s\big\}\,,
\end{split}
\end{equation*}
for some constants $\alpha>0$ and $C_2>1$.   This class is rich; for example it contains the functions that have at least $\alpha-1$ square integrable derivatives. Characterizing error distributions through their Fourier transforms is standard in nonparametric measurement error problems because it is the key to deriving precise asymptotic properties of the estimators.

We consider the sinc kernel $\mathcal{K}$  introduced below \eqref{eq:density}. Using this kernel simplifies our presentation of the theoretical derivations  from two aspects: its Fourier transform simplifies calculations, and it is a so-called infinite order kernel, which implies that the bias of the resulting  nonparametric curve estimators depends only on the smoothness of the target curve.
Thus the sinc kernel has adaptive properties and automatically ensures optimal convergence rates. In contrast, the bias of estimators based on finite order kernels, such as the Gaussian kernel, depends on both the order of the kernel and the smoothness of the target curve, which implies that various smoothness subcases need to be considered when deriving properties.

For any two square integrable functions $f$ and $g$, let $\|f-g\|_2=(\int|f-g|^2)^{1/2}$.
Proposition~\ref{pn:2pre} gives the convergence rate of $\hat{f}_U$, defined at (\ref{eq:density}), to the true density function $f_U$.
\begin{proposition}\label{pn:2pre}
Let $(X_t)_{t\in[0,T]}\in\mathcal{X}(C_1)$ and assume that the errors $\{U_{t_j}\}_{j=0}^n$ satisfy Assumptions~{\rm\ref{as:inde}} and {\rm\ref{as:symm}}, and that  $f_U\in \mathcal{F}_1(\alpha,C_2)$.
Let $\mathcal {P}_1(\alpha,C_1,C_2)$ denote the collection of  models for $(Y_t)_{t\in[0,T]}$ such that  $Y_t=X_t+U_t$.
Under
  Assumption~{\rm\ref{as:space}}, and with the sinc kernel $\mathcal{K}$,
  if $\alpha>3/2$, then for some uniform constant $C>0$ only depending on $\alpha$, $C_1$ and $C_2$,
\[
\sup_{\mathcal{P}_1(\alpha,C_1,C_2)}E\big(\|\hat{f}_U-f_U\|_2^2\big) \leq C(n^{-1}\xi^{-1/2}h^{-1}+n^{-1/2}+\xi+h^{2\alpha-1})\,.\]
\end{proposition}
Proposition~\ref{pn:2pre} shows that the $L_2$ convergence rate of $\hat f_U$ to $f_U$ is affected by $\xi$, the length of each block $S_j$, and the bandwidth $h$.
The next theorem  shows that for appropriate choices of $\xi$ and $h$,  the convergence rate achieves $n^{-1/2}$.

\begin{theorem}\label{tm:1}
If the conditions of Proposition {\rm\ref{pn:2pre}} hold, and we take $\xi\asymp n^{-\delta_1}$ and $h\asymp n^{-\delta_2}$, where $\delta_1>0$ and $\delta_2>0$ are such that $\delta_1+2\delta_2\leq1$, $\delta_1\geq{1}/{2}$ and $\delta_2\geq (4\alpha-2)^{-1}$, then for some uniform constant $C>0$ only depending on $\alpha$, $C_1$ and $C_2$,
\[
\sup_{\mathcal{P}_1(\alpha,C_1,C_2)}E\big(\|\hat{f}_U-f_U\|_2^2\big)\leq Cn^{-1/2}. \]
\end{theorem}

We learn from Theorem~\ref{tm:1} that as long as $\alpha>3/2$, the convergence rate of $\hat f_U$ does not depend on $\alpha$. This is strikingly different from standard nonparametric density estimation problems, where convergence rates typically depend on the smoothness of the target density: the smoother the density, the faster the rates. For example, if we had access to the $U_{t_\ell}-U_{t_j}$'s directly  instead of just $Y_{t_\ell}-Y_{t_j}=X_{t_\ell}-X_{t_j}+U_{t_\ell}-U_{t_j}$, then we could apply the technique suggested by \cite{Meister_2007} and the convergence rate would increase with $\alpha$. However, in our case, the nuisance contribution by the $X_{t_\ell}-X_{t_j}$'s makes it impossible to reach rates faster than $n^{-1/2}$, even if $\alpha$ is very large.  This is demonstrated in the next theorem, which proves that the $n^{-1/2}$ rate derived in Theorem~\ref{tm:1} is minimax optimal.

\begin{theorem}\label{tm:2minmax}
Denote by $\breve{\mathcal{F}}$ the class of all measurable functionals of the data.
Under the conditions in Proposition~{\rm\ref{pn:2pre}},   for some uniform constant $C>0$ only depending on $\alpha$, $C_1$ and $C_2$,
\[
\inf_{\hat{f}\in\breve{\mathcal{F}}}\sup_{\mathcal{P}_1(\alpha,C_1,C_2)}E\big(\|\hat{f}-f_U\|_2^2\big)\geq C n^{-1/2}.
\]
\end{theorem}

\subsection{Estimating the moments of the microstructure noise}\label{se:mome}
We can deduce estimators of the moments of the microstructure noise $U_{t_j}$ from the density estimator derived in Section~\ref{sec:methfU},
but proceeding that way is unnecessarily complex.
Recall from Section~\ref{sec:methfU} that, when $t_\ell\neq t_j$ are close, $Y_{t_\ell}-Y_{t_j}$ behaves approximately like $U_{t_\ell}-U_{t_j}\sim f_{\tilde U}=f_U*f_U$, where $U$ and $\tilde U$ denote generic random variables with the same distribution as, respectively, $U_{t_j}$ and $U_{t_\ell}-U_{t_j}$.
This suggests that we could estimate the moments of $\tilde U$ by the empirical moments of $Y_{t_\ell}-Y_{t_j}$, and, from these, deduce estimators of the moments of $U$.

For each integer $k\geq 1$, let  $M_{U,k}$ and $M_{\tilde{U},k}$ denote the $k$th moment of $U$ and $\tilde U$, respectively. Since $f_U$ is symmetric,  $M_{U,2k-1}$  and $M_{\tilde{U},2k-1}$ are equal to zero for all $k\geq1$, and we only need to estimate even order moments. For each $k\geq1$, we start by estimating $M_{\tilde{U},2k}$ by
\begin{equation*}\label{eq:mtu}
\hat{M}_{\tilde{U},2k}(\xi)= \frac{1}{N(\xi)}\sum_{j=0}^n\sum_{t_\ell\in S_j}(Y_{t_\ell}-Y_{t_j})^{2k}\,.
\end{equation*}
This is directly connected to our frequency domain analysis: it is easily proved that $\hat{M}_{\tilde{U},2k}(\xi)=(-i)^{2k}\{\hat{f}_{\tilde{U}}^{\ft}(0;\xi)\}^{(2k)}$, where $\hat{f}_{\tilde{U}}^{\ft}=(\hat{f}_{U,1}^{\ft})^2$ is an estimator of $f_{\tilde U}^\ft$, with $\hat{f}_{U,1}^{\ft}$ at (\ref{eq:char0}).

Exploiting the fact that $U_{t_\ell}-U_{t_j}\sim f_{\tilde U}$, we can write
$M_{\tilde{U},2k}=\sum_{j=0}^kC_{2k}^{2j}M_{U,2j}M_{U,2k-2j}$,
where $C_{2k}^{2j}=(2k)!/\{(2j)!(2k-2j)!\}$. It can be deduced from there that
\begin{equation*}\label{eq:ture}
M_{U,2k}=\frac{1}{2}\bigg(M_{\tilde{U},2k}-\sum_{j=1}^{k-1}C_{2k}^{2j}M_{U,2j}M_{U,2k-2j}\bigg)\,.
\end{equation*}
Therefore, we can use an iterative procedure  to estimate the $M_{U,2k}$'s. First, for $k=1$, we take $\hat{M}_{U,2}(\xi)=\hat{M}_{\tilde{U},2}(\xi)/2$. Then, for $k>1$, given $\hat{M}_{U,2}(\xi),\ldots, \hat{M}_{U,2(k-1)}(\xi)$ we take
\begin{equation}\label{eq:sample}
\hat{M}_{U,2k}(\xi)=\frac{1}{2}\bigg\{\hat{M}_{\tilde{U},2k}(\xi)-\sum_{j=1}^{k-1}C_{2k}^{2j}\hat{M}_{U,2j}(\xi)\hat{M}_{U,2k-2j}(\xi)\bigg\}\,.
\end{equation}

\begin{remark}
When $k=1$, $M_{U,2}=M_{\tilde{U},2}/2$ is equal to the variance of $U_t$, and our estimator is very similar to the so-called difference-based variance estimator often employed in related nonparametric regression problems; see for example \cite{Buckley} and \cite{HallKayTitt}.
\end{remark}

The next theorem establishes the convergence rate of $\hat{M}_{U,2k}(\xi)$. Its proof follows from  the convergence rates of the $\hat{M}_{\tilde{U},2l}(\xi)$'s. %This result and other results in this and the next sections are derived in the Supplementary Material.
\begin{theorem}\label{tm:3}
Under Assumptions {\rm\ref{as:inde}}--{\rm\ref{as:symm}}, for any integer $k\geq1$, if $E(\int_0^T\mu_s^{2k}ds)<\infty$ and $E(\int_0^T\sigma_s^{2k}ds)<\infty$ hold, and if there exists $p\in(2,3]$ such that $M_{U,2kp}<\infty$, then $ \hat{M}_{U,2k}(\xi)=M_{U,2k}+O_p(n^{-1/2}) $ provided that $\xi=o(n^{-p/\{2(p-1)\}})$.
\end{theorem}

Next we derive the asymptotic joint distribution of the proposed moment estimators. Let  $W=(W_1,\ldots,W_k)^\T$ be a random vector with a $N(0,\Sigma_k)$ distribution, where the $(l_1,l_2)$th element, $l_1, l_2=1,\ldots,k$, of $\Sigma_k$  is equal to
\[
e_{l_1l_2}=\lim_{n\rightarrow\infty}E\bigg(\bigg[\sum_{j=0}^n\sum_{t_\ell\in S_j}\frac{(U_{t_\ell}-U_{t_j})^{2l_1}-M_{\tilde{U},2l_1}}{\{N(\xi)V_{2l_1}(\xi)\}^{1/2}}\bigg]\bigg[\sum_{j=0}^n\sum_{t_\ell\in S_j}\frac{(U_{t_\ell}-U_{t_j})^{2l_2}-M_{\tilde{U},2l_2}}{\{N(\xi)V_{2l_2}(\xi)\}^{1/2}}\bigg]\bigg)\,,
\]
with $V_{2l}(\xi)=\var[ \{N(\xi)\}^{-1/2}\sum_{j=0}^n\sum_{t_\ell\in S_j}\{(U_{t_{\ell}}-U_{t_j})^{2l}-M_{\tilde{U},2l}\}]$.
Recalling that  $N(\xi)\asymp n^2\xi$ and noting that $V_{2l}(\xi)\asymp n\xi$, let
\begin{equation}\label{eq:al}
a_l=\lim_{n\to\infty} \{nV_{2l}(\xi)/N(\xi)\}^{1/2}/2~~~~~(l=1,\ldots,k) \,.
\end{equation}
The next theorem establishes the asymptotic joint distribution of  our moment estimators. It can be used to derive confidence regions for  $(M_{U,2},\ldots,M_{U,2k})^\T$.
\begin{theorem}\label{tm:limitdist}
Under Assumptions {\rm\ref{as:inde}}--{\rm\ref{as:symm}}, for any integer $k\geq1$, if $E(\int_0^T\mu_s^{2k}ds)<\infty$ and $E(\int_0^T\sigma_s^{2k}ds)<\infty$, and if there exists $p\in(2,3]$ such that $M_{U,2kp}<\infty$, then
\[
n^{1/2}\{\hat{M}_{U,2}(\xi)-M_{U,2},\ldots,\hat{M}_{U,2k}(\xi)-M_{U,2k}\}^\T\xrightarrow{d}Q=(Q_{1},\ldots, Q_{k})^\T~~\textrm{as}~~n\rightarrow\infty\,,
\]
provided that $\xi=o(n^{-p/\{2(p-1)\}})$, where
 \begin{equation}\label{eq:u}
         Q_1 = a_1W_1,~~~~~~Q_l = a_lW_l-\sum_{j=1}^{l-1}C_{2l}^{2j}M_{U,2l-2j}Q_j~~~~(2\leq l\leq k)\,.
                          \end{equation}
\end{theorem}

\subsection{Efficient integrated volatility estimation}
\label{s4}
\label{se:volatility}
We have demonstrated that our frequency domain analysis can be used to estimate the error density, which is difficult to estimate in the time domain.
Since frequency domain approaches are unusual in  high-frequency financial data analysis,
a natural question is whether they can lead to an efficient estimator of the integrated volatility $\int_0^T\sigma_t^2~dt$, with $\sigma_t$ as in \eqref{eq:model}. The integrated volatility is a key quantity of interest in high-frequency financial data analysis; it represents the variability of a process over time. It is well known that in cases like ours where the data are observed with microstructure noise, the integrated volatility cannot be estimated using standard procedures, which are dominated by contributions from the noise. There, one way of removing the bias caused by the noise is through multiscale techniques; see \cite{OSP_2019} for a very nice description. \cite{Zhang_2006_Bern} and \cite{TaoWangZhou_2013} have successfully applied those methods to correct the bias of estimators in the time domain, and \cite{OSP_2019} have proposed a consistent discrete frequency domain estimator in the case where the data are observed at equispaced times. Next we show that these techniques can be applied in our continuous frequency domain context too, even if the observation times are  not restricted to be equispaced.

The real part of the empirical characteristic function $n^{-1}\sum_{i=1}^n\exp\{i s(Y_{t_i}-Y_{t_{i-1}})\}$ is such that
\begin{equation}
\begin{split}
\sum_{j=1}^n\cos\{s(Y_{t_j}-Y_{t_{j-1}})\}
=&~\sum_{j=1}^n\cos\{s(U_{t_j}-U_{t_{j-1}})\}-\frac{s^2}{2}f_{\tilde{U}}^\ft(s)\int_0^T\sigma_t^2~dt\\
&+O_p(n^{-1/2})\,.\label{phiDiffY}
\end{split}
\end{equation}
The second term on the right hand side of \eqref{phiDiffY} contains the integrated volatility, but the first term dominates  because its mean is $nf_{\tilde U}^\ft(s)$.
This suggests that the integrated volatility could be estimated from $\sum_{i=1}^n\exp\{i s(Y_{t_i}-Y_{t_{i-1}})\}$,  if we could eliminate that first term.
This can be done by applying, to the frequency domain, the multiscale technique used by \cite{Zhang_2006_Bern} and \cite{TaoWangZhou_2013}. We define a function $G(s)$ which combines the empirical characteristic function calculated at different sampling frequencies, in such a way as to eliminate the first term on the right hand side of \eqref{phiDiffY} while keeping the second.  For $N=\lfloor (n+1)^{1/2} \rfloor$,
we define
\[G(s)=\sum_{m=1}^N a_m \hat\phi^{K_m}(s)+\zeta \{\hat\phi^{K_1}(s)-\hat\phi^{K_2}(s)\}\,,\]
where, as in \cite{Zhang_2006_Bern}, $K_m=m$, $a_m=12K_m(m-N/2-1/2)/\{N(N^2-1)\}$, $\zeta=K_1K_2/\{(n+1)(K_2-K_1)\}$, and where
$\hat\phi^{K_m}(s)=K_m^{-1}\sum_{\ell=K_m}^n \exp\{i s(Y_{t_\ell}-Y_{t_{\ell-K_m}})\}$. We could also select $K_m$, $a_m$ and $\zeta$ as in \cite{TaoWangZhou_2013}.

The following proposition shows that the real part of $G(s)$, $\Rebu\{G(s)\}$, can be used to approximate the second term of \eqref{phiDiffY}.
\begin{proposition}\label{pn:diff}
Under Assumptions {\rm\ref{as:inde}}--{\rm\ref{as:symm}}, if $E(U_t^2)<\infty$, then there exist second-order differentiable
functions $\tau_1(s)$ and $\tau_2(s)$ such that for any $s\in R$
\[
\bigg|\Rebu\{G(s)\}+\frac{s^2}{2}f_{\tilde{U}}^{\ft}(s)\int_0^T\sigma_t^2~dt\bigg|=\tau_1(s)\, O_p(n^{-1/4})+\tau_2(s)\, O_p(n^{-1/2}),
\]
where the terms $O_p(n^{-1/4})$ and $O_p(n^{-1/2})$ are independent of $s$, and where $\lim_{s\to 0}|\tau_1''(s)|\leq C$ and $\lim_{s\to 0}|\tau_2''(s)|\leq C$ for some positive constant $C$.
\end{proposition}

Since the function $G(s)$ depends only on the data, it is completely known. Moreover we have seen in Section~\ref{se:mome} that we could estimate  $f_{\tilde{U}}^{\ft}(s)$ by $\hat{f}_{\tilde{U}}^{\ft}(s;\xi)=\{\hat{f}_{U,1}^{\ft}(s;\xi)\}^2$. Finally, although the proposition holds for all $s\in R$, the remainders are smaller when $f_{\tilde{U}}^{\ft}(s)$ is close to one, especially since $\hat{f}_{\tilde{U}}^{\ft}(s;\xi)$ is more reliable in that case too. Therefore, for $\hat{f}_{\tilde{U}}^{\ft}(s;\xi)$ close to one,
Proposition~{\ref{pn:diff}} can be used to compute an estimator of the integrated volatility, $\int_0^T\sigma_t^2~dt$.
We propose a regression type approach as follows. For some $s_1,\ldots,s_m$ such that $\hat{f}_{\tilde{U}}^{\ft}(s;\xi)$ is close to one we consider the fixed design regression problem
\begin{equation*}\label{eq:regress}
\Rebu\{G(s_j)\}=-\frac{s_j^2}{2}\hat{f}_{\tilde{U}}^{\ft}(s_j;\xi)\cdot \beta+\epsilon_j~~~~~(j=1,\ldots,m)\,,
\end{equation*}
where $\epsilon_j$ represents the regression error and $\beta=\int_0^T\sigma_t^2~dt$.  Applying a linear regression of $\Rebu\{G(s_j)\}$  on $-s_j^2\hat{f}_{\tilde{U}}^{\ft}(s_j;\xi)/2$, we  estimate $\int_0^T\sigma_t^2~dt$ by $\hat \beta$, the least squares estimator of $\beta$.

For any fixed $s\in R$, it can be shown that  $|\hat{f}_{\tilde{U}}^{\ft}(s;\xi)-f_{\tilde{U}}^{\ft}(s)|=O_p(n^{-1/2}+\xi)$. If we select   $\xi=O(n^{-1/4})$ in (\ref{eq:sj}), then Proposition \ref{pn:diff} still holds if we replace $f_{\tilde{U}}^{\ft}(s)$ by  $\hat{f}_{\tilde{U}}^{\ft}(s;\xi)$. The next result establishes the convergence rate of $\hat\beta$.
\begin{theorem}\label{cor1}
Under Assumptions {\rm\ref{as:inde}}--{\rm\ref{as:symm}}, if $E(U_t^2)<\infty$ and selecting $\xi=O(n^{-1/4})$, then
\[
\hat{\beta}-\int_0^T\sigma_t^2~dt=O_p(n^{-1/4})\,.
\]
\end{theorem}
\noindent
The convergence rate stated in Theorem~\ref{cor1} is optimal in the sense of \cite{GloterJacod_2001}, and is the same as for the time-domain estimator of \cite{TaoWangZhou_2013}.
Hence, our frequency domain method is rate-efficient for estimating the integrated volatility.

\section{Numerical study}
\label{s3}
\subsection{Practical implementation of the density estimator}\label{ParamSelec}
To compute the density estimator $\hat f_U$ at \eqref{eq:density}, we need to choose the bandwidth $h$ and the parameter $\xi$.
In our problem, doing this is much more complex than that for standard nonparametric density estimators, since, unlike there, we  do not have direct access to  data from our target $f_U$.
Therefore, we cannot use existing smoothing parameter selection methods, which all require noise free data.
Moreover, unlike in standard nonparametric problems, we do not have access to a formula measuring the distance between $f_U$ and its estimator.

Similar difficulties arise in the classical errors-in-variables problem, where one is interested in a density $f_V$, but  observes only data on $W=V+\varepsilon$, where $V\sim f_V$ is independent of $\varepsilon\sim f_\varepsilon$ with $f_\varepsilon$ known.  \cite{DelaigleHall_2008} proposed to choose  the bandwidth $h$ using a method called simulation approximation. Instead of computing $h$ for $f_V$ they approximate $h$ by extrapolating the bandwidths for estimating two other densities, $f_{1}$ and $f_{2}$, that are related to $f_V$. The rationale of the extrapolation scheme is to exploit the analogous  relationships between $f_1$, $f_2$,  and $f_V$.

Our problem is different because our estimator of $f_U$ is completely different from that in \cite{DelaigleHall_2008}, and we do not know $f_X$.  Therefore, we cannot apply their method directly,
and here we propose a method with the same spirit.
In particular, we consider two density functions  $f_{1}(\cdot)=2f_U(\surd2\, \cdot)*f_U(\surd2\, \cdot)$ and $f_{2}=2f_{1}(\surd2\, \cdot)*f_{1}(\surd2\,
\cdot)$, and remark that the way in which $f_2$ and $f_1$ are connected mimics the way in which $f_1$ and $f_U$ are connected.
As shown in the Appendix,   data from  both $f_1$ and $f_2$ can be made accessible, so that one may perform bandwidth selection for estimating them.
We then choose the bandwidth for estimating $f_U$ by an  extrapolation with a ratio adjustment from  the bandwidths for estimating $f_1$ and $f_2$.

The algorithms for the bandwidth selection are given in the Appendix, and
Algorithm~\ref{alg:band} summarizes the main steps.
Specifically,
for $k=1,2$, our procedure requires the construction of variables $\Delta Y_{j,\ell}^{_k}$ and times points $t_j^{_k}$, which are defined in, respectively, Algorithms~\ref{alg:DeltaY} and \ref{alg:DeltaY2} in the Appendix. Step (c) of Algorithm~\ref{alg:band} requires choosing values of   $(h,\xi)$, say $(h_{k},\xi_{k})$, for $k=1,2$,
for estimating  $f_{k}$ by $\hat f_{k}$, where $\hat f_{k}$ denotes our density estimator at \eqref{eq:density} applied to the $\Delta Y_{j,\ell}^{_k}$'s.

The idea is that if we knew $f_{1}$, we would choose $(h_{1},\xi_{1})$ by minimizing the integrated squared error of $\hat f_{1}$, i.e.
$
(h_{1},\xi_{1})={\rm argmin}_{(h,\xi)} \int \{ \hat f_{1}(x;\xi) - f_{1}(x)\}^2\,dx.%\label{hxi1}
$
In practice we do not know $f_{1}$, but we can construct a relatively accurate estimator  of it,
the standard kernel density estimator $\tilde f_{1}$ of $f_{1}$ applied  to the data $\Delta_{Y,j}=(Y_{t_{j+1}}-Y_{t_{j}})/\surd 2\approx (U_{t_{j+1}}-U_{t_{j}})/\surd 2$  $(j=0,\ldots,n-1)$. This is computed by using the Gaussian kernel with bandwidth selected by the method of \cite{SheatherJones}.
Using arguments similar to those in \cite{Delaigle_2008}, under mild conditions,
$\tilde f_{1}(x)=f_{1}(x)+O_p(n^{-2/5})$, whereas, with the best possible choice of $(h,\xi)$,  $\hat f_{1}(x)=f_{1}(x)+O_p(n^{-1/4})$, where the rate $n^{-1/4}$ cannot be improved. Thus,  $\tilde f_{1}$ converges to $f_{1}$ faster than   $\hat f_{1}$  does.  This motivates us to approximate $(h_{1},\xi_{1})$ defined above by
\begin{equation}
(h_{1},\xi_{1})=\arg\min_{(h,\xi)} \int \{ \hat f_{1}(x;\xi) -\tilde f_{1}(x)\}^2\,dx\,.\label{hathxi}
\end{equation}

Paralleling the arguments in \cite{DelaigleHall_2008}, it is more important to extrapolate the bandwidth $h$ than $\xi$. Motivated by their results, we take $\xi_{2}=\xi_{1}$. To choose $h_{2}$, let $\tilde f_{2}$ be the standard kernel density estimator with the Gaussian kernel and bandwidth selected by the method of \cite{SheatherJones}, applied to the data $\Delta_{Y,j,2}=(\Delta_{Y,j}-\Delta_{Y,k(j)})/\surd 2$, where $k(j)$ is chosen at random from $0,\ldots,n-1$. We choose $\Delta_{Y,j,2}$ in this way rather than $\Delta_{Y,j,2}=(\Delta_{Y,j}-\Delta_{Y,j+2})/\surd 2$, to prevent accumulated residual $X_t$ effects.
Then we take
\begin{equation}
h_{2}=\arg\min_{h}   \int \{ \hat f_{2}(x;\xi_{1}) - \tilde f_{2}(x)\}^2\,dx\,. \label{hath2}
\end{equation}
Since $\tilde f_1$ and $\tilde f_2$ converge faster than $\hat f_1$ and $\hat f_2$, they can be computed using less data than the latter. Therefore, when the time points are widely unequispaced, for computing the $\tilde f_k$'s we suggest using only a fraction, say one quarter, of the $\Delta_{Y_j}$'s corresponding to the smallest $t_j-t_{j+1}$'s. That is, we use less but more accurate data for computing the $\tilde f_k$'s.

Finally, as described in step (d) of Algorithm~\ref{alg:band}, we obtain our bandwidth for estimating $f_U$ by an extrapolation with a ratio adjustment: $\hat h = h_{1}^2/ h_{2}$, and we take $\hat\xi=\xi_1$.
This method is not guaranteed to give the best possible bandwidth  for estimating $f_U$, but is a sensible approximation for a problem that seems otherwise very hard, if not impossible, to solve.  Theorem~\ref{tm:1} implies that we have a lot of flexibility for choosing $h$, but it is impossible to know if our bandwidth lies in the optimal range without knowing the exact order of $h_1$ and $h_2$. However, we cannot determine these  orders without deriving complex second-order asymptotic results.

\subsection{Simulations}\label{sec:simul}

\begin{figure}[t]
\begin{center}
\vspace*{-1cm}
\makeatletter\def\@captype{figure}\makeatother
\centering
\hspace*{-.2cm}
\includegraphics[width=5.2cm]{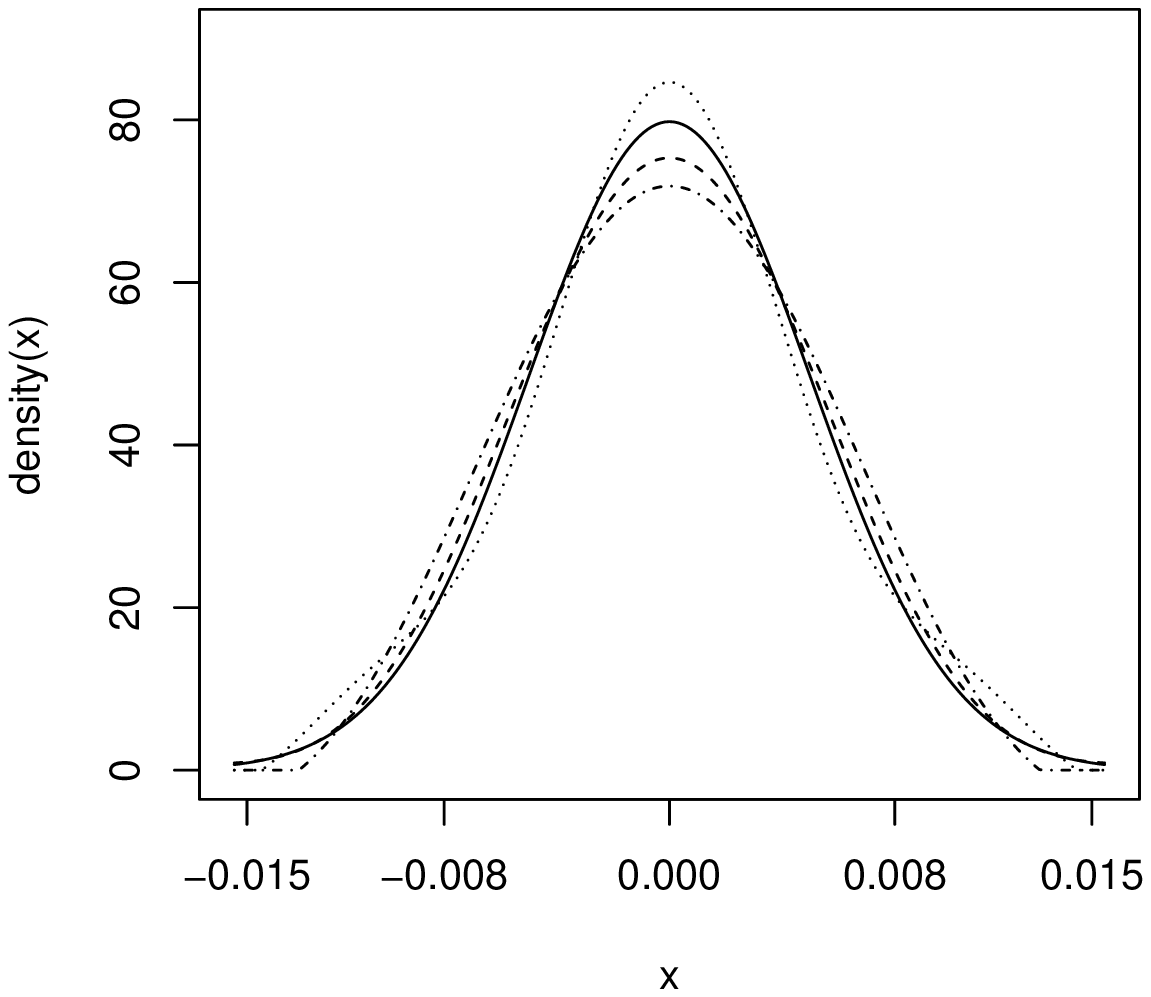}\hspace*{-.4cm}
\includegraphics[width=5.2cm]{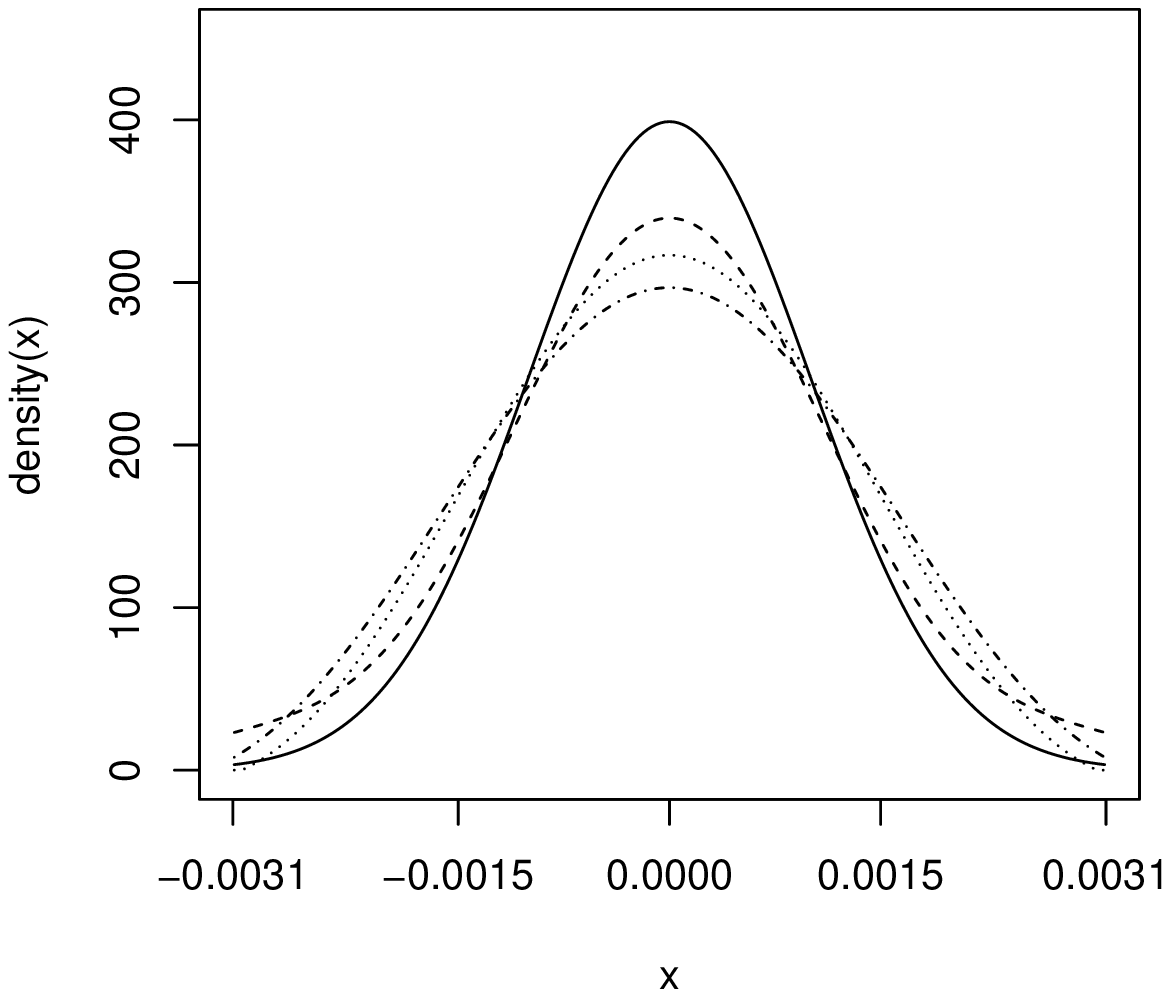}\hspace*{-.4cm}
\includegraphics[width=5.2cm]{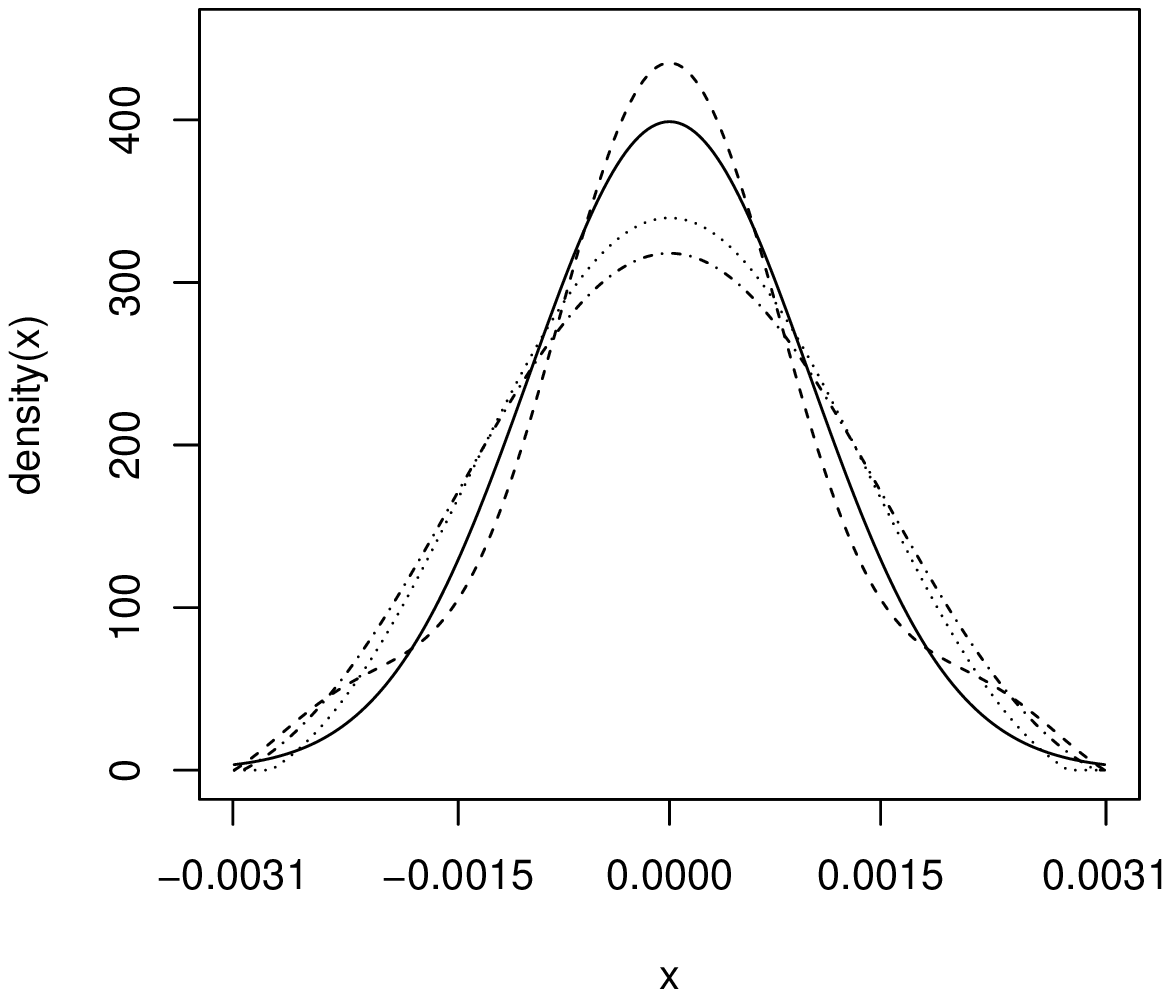}\\[-1cm]
\hspace*{-.2cm}
\includegraphics[width=5.2cm]{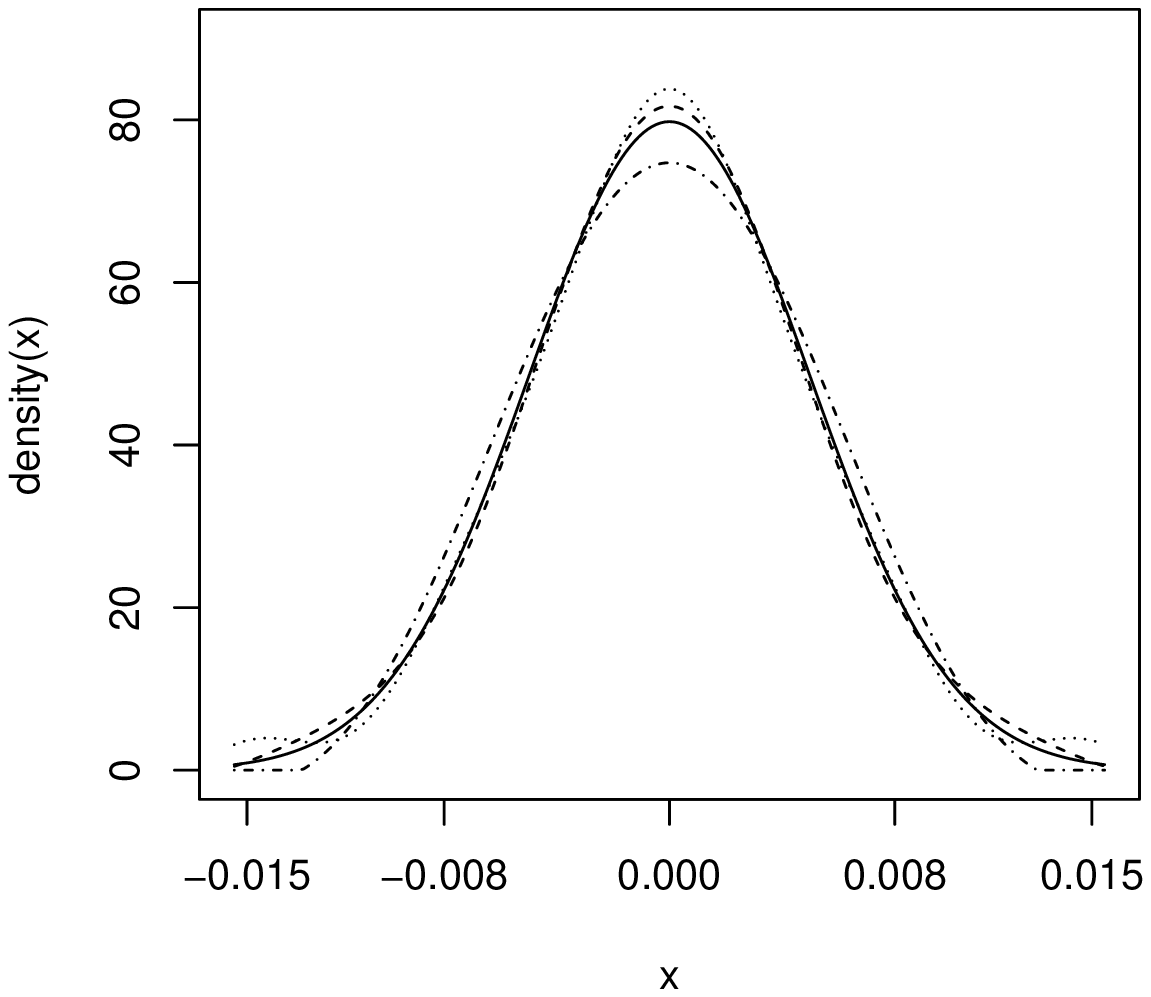}\hspace*{-.4cm}
\includegraphics[width=5.2cm]{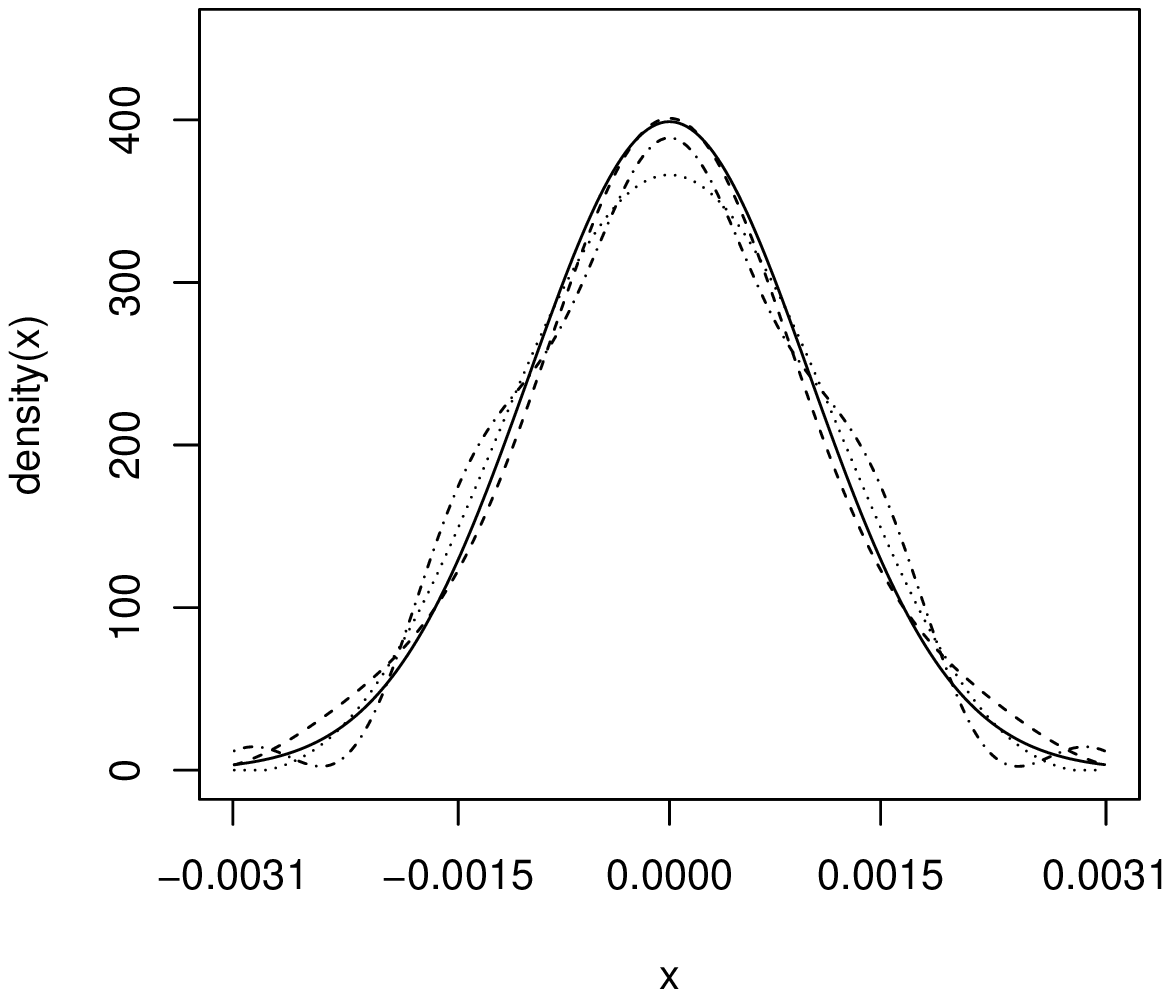}\hspace*{-.4cm}
\includegraphics[width=5.2cm]{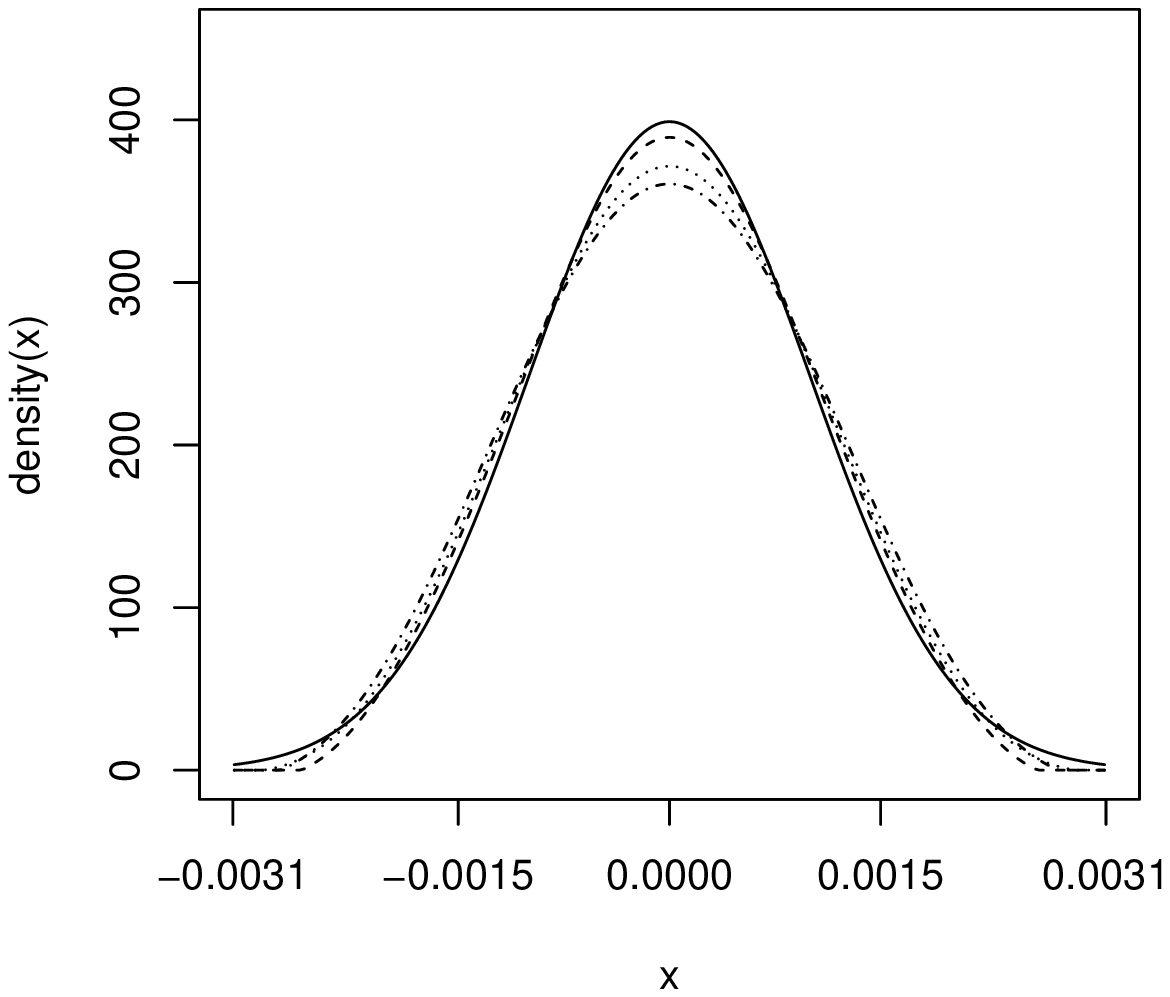}\\[-.6cm]
\caption{Estimator $\hat f_U(x)$ at \eqref{eq:density} in the case of normal errors, for three samples corresponding to the  first, - - -, second, $\cdots$, and third, -$\cdot$-, quartiles of the integrated square errors of estimators computed from data from model (i) with $\sigma_U=0$$\cdot$$005$, colomun 1,  and with $\sigma_U=0$$\cdot$$001$, column 2, and model (ii) with $\sigma_U=0$$\cdot$$001$, column 3,  when  $\Delta s=30$, row 1, and $\Delta s=5$, row 2.  The continuous line depicts the true~$f_U(x)$.}\label{F:Fig1}
\end{center}
\end{figure}

We applied our methodology to  data simulated from stochastic volatility models. We generated the data $Y_{t_0},\ldots,Y_{t_n}$ as in \eqref{eq:finalmodel}. Taking the convention that a financial year has $252$ active days, we took $t\in[0,T]$ with $T=1/252$; which represents one day of financial activity. We took time points every $\Delta s$ seconds, where $\Delta s=30$, $5$ and $1$. Taking the convention of $6.5$ business hours for a trading day,  this means that we  took the $t_j$'s to be equally spaced by $\Delta s/(252\times60\times60\times6.5)$, and that $n$ was equal to $60\times60\times6.5/\Delta s$. We generated the microstructure noise $U_t$ according to a normal  or a scaled $t$-distribution, and for the $X_t$'s, we used the  Heston model
\[
dX_t=\sigma_t~dB_t\,,\ \
d\sigma_t^2=~\kappa(\tau-\sigma_t^2)~dt+\gamma \sigma_t~dW_t\,,
\]
where $E(dB_t\, dW_t)=\rho\, dt$ and $\kappa$, $\tau$, $\gamma$ and $\rho$ are parameters.
Like \cite{AitYu_2009_AOAS}, we set the drift part of $X_t$ to zero. The impact due to the drift function is asymptotically negligible; see, for example, \cite{Xiu_2010_JOE}.
We used two models, with values similar to those used by \cite{AitYu_2009_AOAS}, which reflect practical scenarios in finance; see also \cite{Xiu_2010_JOE} and \cite{LiuTang_2014_JOE}:
(i) $(\kappa, \tau, \gamma, \rho)=(6,0.16,0.5,-0.6)$ and
(ii) $(\kappa, \tau, \gamma, \rho)=(4, 0.09, 0.3, -0.75)$.
In each case we took $X_0=\log(100)$ and considered $U_t\sim N(0,\sigma_U^2)$ and $U_t\sim \sigma_U\, t(8)$, where $\sigma_U=0.001$ and $\sigma_U=0.005$.
Typical  $X_t$ and $Y_t$ paths for each model plotted in the Supplementary Material show that the $X_t$ paths are with smaller variation in model (ii) than in model (i). The $X_t$ paths with smaller variation have less nuisance impact on estimators of $U_t$-related quantities. Thus, estimating the moments and density of $U_t$ should be easier in model (ii).

In each setting,  we generated $1000$ samples of the form $Y_{t_0},\ldots,Y_{t_n}$, and applied our estimator of the density $f_U$ to each sample, obtaining in this way 1000 density estimators $\hat f_U$ computed as in \eqref{eq:density}. We chose the smoothing parameters as in Section~\ref{ParamSelec}, and took the sinc kernel defined below \eqref{eq:density}. However, while this kernel guarantees optimal theoretical properties, in practice it produces negative wiggles in the tails, which we truncate to zero since $f_U$ is a density. In the Supplementary Material, we show the results obtained when using the Gaussian kernel, which suggest that, overall, the sinc kernel works better, but the Gaussian kernel produces more attractive estimators in the tails. In cases where the sample has ties, for example when the sample size is large and the data are observed only with a few significant digits, like in our real data example, the wiggles of the sinc kernel make it perform too poorly and significantly better results can be obtained when using the Gaussian kernel; see Section~\ref{sec:realdata}.

For each estimator, we computed the integrated squared error, $\int (\hat f_U-f_U)^2$, and their median, first and third quartiles are reported in Table~\ref{table:Simul}.  In Fig.~\ref{F:Fig1} we depict, for selected settings with normal error, the estimated curves $\hat f_U$ computed from the samples corresponding to those three quartiles. Our results indicate that our density estimator works well. For a given setting, error densities with larger variances are easier to estimate. Fig.~\ref{F:Fig1} shows that our estimator improves as the sample size increases, that is as $\Delta s$ decreases. The estimated densities are better in model (ii) than in model (i), as expected. While it is difficult to compare estimators of different target densities, the figures suggest that the difficulty in estimating the error densities depends more on the smoothness of $X_t$ than on the error type. This reflects the fact that the smoothness of the error density has no first order impact on the quality of estimators, as indicated by Theorems~\ref{tm:1} and \ref{tm:2minmax}.

We also applied the moment estimators from Section~\ref{se:mome}  to data simulated from a rescaled version of our two models; the rescaling was applied to avoid working with numbers all numerically rounded to zero. Namely, we replaced  $(\tau, \gamma, \sigma^2_U,X_0)$ in our models by $(c^2\tau, c\gamma, c^2 \sigma^2_U,cX_0)$, where $c=100$.  We present the results in Table~\ref{table:SimulMom}. The results indicate that, as expected, performance improves as the sample size increases. Moreover, the performance is best for lower order moments, for higher noise levels, and for model (ii). For a given error variance, the moments are easier to recover when the errors have a rescaled Student distribution than when they have a normal distribution.

\begin{table}[t]
\begin{center}
 \caption{Median integrated squared error $[$first quartile, third quartile$]$ of $\hat f_U$ at \eqref{eq:density}, calculated from $1000$ simulated samples from models (i) and (ii)}\label{table:Simul}
\begin{tabular}{llllll}
%\hline
%\hline
&&\multicolumn{2}{c}{Normal errors}&\multicolumn{2}{c}{Student errors}\\
$\Delta s$&$\sigma_U$&Model (i) & Model (ii) & Model (i) & Model (ii) \\
%\hline
$30$&0$\cdot$005&0$\cdot$29 [0$\cdot$13,0$\cdot$59]&0$\cdot$26 [0$\cdot$10,0$\cdot$58]&0$\cdot$37 [0$\cdot$19,0$\cdot$79]&0$\cdot$36 [0$\cdot$16,0$\cdot$76]\\
&0$\cdot$001&9$\cdot$98 [5$\cdot$39,16$\cdot$6]&5$\cdot$63 [2$\cdot$90,10$\cdot$2]&9$\cdot$82 [5$\cdot$8,15$\cdot$59]&6$\cdot$14 [3$\cdot$28,10$\cdot$4]\\
$5$ &0$\cdot$005&0$\cdot$10 [0$\cdot$04,0$\cdot$27]&0$\cdot$10 [0$\cdot$04,0$\cdot$29]&0$\cdot$12 [0$\cdot$06,0$\cdot$27]&0$\cdot$14 [0$\cdot$07,0$\cdot$33]\\
 &0$\cdot$001&1$\cdot$46 [0$\cdot$64,3$\cdot$21]&1$\cdot$02 [0$\cdot$41,2$\cdot$18]&1$\cdot$59 [0$\cdot$78,2$\cdot$97]&1$\cdot$16 [0$\cdot$56,2$\cdot$19]\\
$1$ &0$\cdot$005&0$\cdot$04 [0$\cdot$01,0$\cdot$14]&0$\cdot$04 [0$\cdot$01,0$\cdot$15]&0$\cdot$05 [0$\cdot$02,0$\cdot$14]&0$\cdot$05 [0$\cdot$02,0$\cdot$13]\\
 &0$\cdot$001&0$\cdot$32 [0$\cdot$12,1$\cdot$01]&0$\cdot$25 [0$\cdot$10,0$\cdot$82]&0$\cdot$42 [0$\cdot$20,0$\cdot$95]&0$\cdot$36 [0$\cdot$15,0$\cdot$86]\\
%\hline
\end{tabular}
\end{center}
\end{table}

Finally we applied our volatility estimator in each of our four settings, and compared it with the estimator of \cite{Zhang_2006_Bern}.   For our method we took $\xi=t_2-t_1$ and to choose $s_1,\ldots,s_m$,  we took $m=50$ equispaced points located between $0$ and $S$, where $S$ was the largest number for which $\hat{f}_{\tilde{U}}^{\ft}(S;\xi)\geq 0.99$.
The results are presented in Table~\ref{table:SimulVolAlt}. In the Supplementary Material, we also show the first three quartiles of  $10^2\times$ the relative absolute deviation, $|\hat \beta-\beta|/\beta$, of both estimators. Together, the results indicate that overall both estimators gave similar results. In the case where the error variance was large, our estimator tended to work a little better, and in the case where the error variance was small, the estimator from \cite{Zhang_2006_Bern} worked a little better. Specifically, there it was a little less biased, but ours was a little less variable. As expected, all estimators improved as the sample size increased.
The estimators worked a little better when the error variance was small: the estimators were a little more biased but significantly less variable. The integrated volatility is an $X_t$-related quantity, and is therefore easier to estimate when there are less errors contaminating the $X_t$'s. Depending on the case, either the Student errors or the normal errors gave better results.

\begin{table}[t]
\begin{center}
\caption{$10^2\times$ bias (standard deviation) of the 1000 relative deviations $(\hat{M}_{U,{2k}}-M_{U,{2k}})/M_{U,{2k}}$ of our estimator $\hat{M}_{U,2k}$ at \eqref{eq:sample}  calculated from $1000$ simulated samples from models (i) and (ii) and computed using $\xi=t_2-t_1$, for $k=1$ and $2$}\label{table:SimulMom}
\setlength\tabcolsep{4.5pt}
\begin{tabular}{lcccrrrrrrrr}
&&&  \multicolumn{2}{c}{$\Delta s$=30}& \multicolumn{2}{c}{$\Delta s$=5}& \multicolumn{2}{c}{$\Delta s$=1}\\
Errors & Model&$\sigma_U$&$k=1$ & $k=2$& $k=1$ & $k=2$ & $k=1$ & $k=2$\\%& $k=3$\\

Normal& (i) &$0$$\cdot$$005$&
1$\cdot$7 (6$\cdot$4)&3 (18)&0$\cdot$27 (2$\cdot$5)&0$\cdot$53 (7)&0$\cdot$08 (1$\cdot$1)&0$\cdot$15 (3$\cdot$3)\\

&&$0$$\cdot$$001$&
41 (15)&99 (51)&6$\cdot$8 (3$\cdot$6)&14 (9$\cdot$6)&1$\cdot$4 (1$\cdot$2)&2$\cdot$8 (3$\cdot$5)\\

&(ii)&$0$$\cdot$$005$&
1 (6$\cdot$3)&1$\cdot$6 (18)&0$\cdot$15 (2$\cdot$5)&0$\cdot$29 (7)&0$\cdot$06 (1$\cdot$1)&0$\cdot$11 (3$\cdot$3)\\

&&$0$$\cdot$$001$&
23 (10)&51 (31)&3$\cdot$8 (2$\cdot$9)&7$\cdot$8 (8$\cdot$1)&0$\cdot$79 (1$\cdot$2)&1$\cdot$6 (3$\cdot$4)\\

Student &  (i)&$0$$\cdot$$005$&
1$\cdot$4 (7$\cdot$7)&0$\cdot$5 (37)&0$\cdot$28 (3)&0$\cdot$13 (17)&0$\cdot$03 (1$\cdot$4)&-0$\cdot$2 (9$\cdot$9)\\

&&$0$$\cdot$$001$&

30 (14)&45 (48)&5$\cdot$2 (3$\cdot$6)&6$\cdot$9 (18)&1 (1$\cdot$4)&1$\cdot$1 (9$\cdot$9)\\

& (ii)&$0$$\cdot$$005$&

0$\cdot$9 (7$\cdot$7)&-0$\cdot$2 (37)&0$\cdot$2 (3)&0$\cdot$01 (17)&0$\cdot$01 (1$\cdot$4)&-0$\cdot$23 (9$\cdot$9)\\

&&$0$$\cdot$$001$&
17 (10)&23 (43)&3 (3$\cdot$2)&3$\cdot$8 (17)&0$\cdot$57 (1$\cdot$4)&0$\cdot$49 (9$\cdot$9)\\

\end{tabular}
\end{center}
\end{table}

\begin{table}[t]
\begin{center}
\caption{$10^2\times$ bias (standard deviation) of $(\hat \beta-\beta)/\beta$ using our estimator of $\int_0^T \sigma_t^2$\,dt from Section~\ref{se:volatility}, denoted by Ours in the table, and that of \cite{Zhang_2006_Bern}, denoted by Zhang in the table, calculated from $1000$ simulated samples from models (i) and (ii)}\label{table:SimulVolAlt}
\begin{tabular}{lllrrrr}
&&&\multicolumn{2}{c}{Normal errors}&\multicolumn{2}{c}{Student errors}\\
$\sigma_U$&$\Delta s$&$\hat\beta$&Model (i)&Model (ii)&Model (i)&Model (ii)\\
$0$$\cdot$$005$&30&Ours&-1$\cdot$12 (25$\cdot$9)&0$\cdot$42 (29$\cdot$6)&-1$\cdot$50 (28$\cdot$0)&0$\cdot$71 (33$\cdot$6)\\
&&Zhang&-2$\cdot$87 (26$\cdot$0)&-2$\cdot$91 (29$\cdot$7)&-3$\cdot$89 (28$\cdot$2)&-3$\cdot$74 (33$\cdot$7)\\
&5&Ours&-1$\cdot$13 (17$\cdot$5)&-0$\cdot$65 (20$\cdot$1)&-0$\cdot$78 (18$\cdot$3)&-0$\cdot$53 (21$\cdot$8)\\
&&Zhang&-1$\cdot$38 (17$\cdot$5)&-1$\cdot$18 (20$\cdot$1)&-1$\cdot$15 (18$\cdot$3)&-1$\cdot$25 (21$\cdot$8)\\
&1&Ours&-0$\cdot$61 (11$\cdot$2)&-0$\cdot$78 (12$\cdot$7)&-0$\cdot$59 (12$\cdot$1)&-0$\cdot$59 (14$\cdot$4)\\
&&Zhang&-0$\cdot$64 (11$\cdot$2)&-0$\cdot$87 (12$\cdot$7)&-0$\cdot$65 (12$\cdot$1)&-0$\cdot$73 (14$\cdot$5)\\
$0$$\cdot$$001$&30&Ours&-5$\cdot$29 (21$\cdot$4)&-4$\cdot$38 (21$\cdot$8)&-6$\cdot$25 (21$\cdot$6)&-5$\cdot$40 (22$\cdot$1)\\
&&Zhang&-2$\cdot$80 (22$\cdot$3)&-2$\cdot$82 (22$\cdot$4)&-4$\cdot$25 (22$\cdot$5)&-4$\cdot$21 (22$\cdot$7)\\
&5&Ours&-3$\cdot$17 (14$\cdot$4)&-2$\cdot$54 (14$\cdot$7)&-2$\cdot$30 (14$\cdot$5)&-1$\cdot$80 (14$\cdot$7)\\
&&Zhang&-1$\cdot$73 (14$\cdot$8)&-1$\cdot$71 (14$\cdot$9)&-1$\cdot$17 (14$\cdot$8)&-1$\cdot$16 (14$\cdot$9)\\
&1&Ours&-1$\cdot$02 (9$\cdot$62)&-0$\cdot$72 (9$\cdot$71)&-1$\cdot$15 (9$\cdot$67)&-0$\cdot$91 (9$\cdot$76)\\
&&Zhang&-0$\cdot$31 (9$\cdot$75)&-0$\cdot$32 (9$\cdot$79)&-0$\cdot$63 (9$\cdot$76)&-0$\cdot$62 (9$\cdot$80)
%\hline
\end{tabular}
\end{center}
\end{table}

\subsection{Real data analysis}\label{sec:realdata}
We applied our procedure for analyzing the high-frequency price data of Microsoft Corporation from March 19 to April 2, 2013, ten trading days, available from the Trade and Quote database. We took the $Y_t$'s equal to the log prices.  Following  \cite{Barndorff 2011}, we pre-processed the data by  deleting entries that have  0 or negative prices, deleting entries with negative Correlation Indicator, deleting entries with  a letter code in COND, except for E or F,  deleting entries outside the period 9:30 a.m. to 4:00 p.m., and using the median price if there were multiple entries at the same time.
In this example, the sample size is very large, which entails a large number of ties among the $Y_t$'s. As a result, the sinc kernel produced too wiggly estimates,  regarding the data as coming from a multimodal density with modes located at the ties. The  oscillation problems of this kernel are well known, but here they are exacerbated by the ties, so we used the standard Gaussian kernel, which is less affected by them. The only small adjustment we had to make was to break the ties by adding a small perturbation $\epsilon_t\sim N(0,a_j^2)$ to the $\Delta_{Y,j}$'s when computing the bandwidth of the standard kernel estimator $\tilde f_{1}(x)$ used in our bandwidth selection procedure, where $2 a_j$ was equal to the maximum of the distance between $\Delta_{Y,j}$ and its first smaller and larger nearest neighbours.

\begin{figure}[t]
\begin{center}
\vspace*{-1.2cm}
\makeatletter\def\@captype{figure}\makeatother
\centering
\hspace*{-.5cm}
\includegraphics[width=0.5\textwidth]{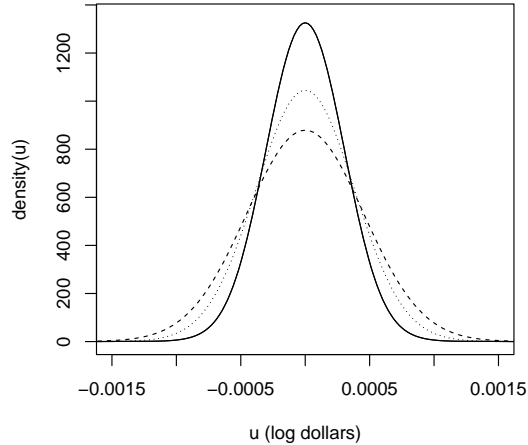}\\[-.6cm]
\caption{Estimated densities  of  the errors contaminating the log prices of the Microsoft Corporation data in March 20 ---, March  28, $\cdots$ and April 2, $- - -$, 2013. }\label{F:Fig3}
\end{center}
\end{figure}

Fig.~\ref{F:Fig3} shows the error densities estimated by our method  for three trading days in 2013:  March 20,  March 28, and April 2.
In this example,  the magnitude of the errors is  about $10^{-3}$ smaller than that of the log prices themselves, but their aggregated impact on quantities such as integrated volatility is substantial.
For example, for those three trading days the  realized volatility was respectively 3$\cdot$5, 4$\cdot$5 and 3$\cdot$5 times $10^{-4}$, which is dominated by contributions from the errors. Indeed, for the same days, our error-corrected estimator was respectively 0$\cdot$5, 0$\cdot$5 and 0$\cdot$3 times $10^{-4}$, and the estimator of \cite{Zhang_2006_Bern} was respectively 0$\cdot$5, 0$\cdot$6 and 0$\cdot$4 times $10^{-4}$.

Interestingly, even for this short period of ten trading days, the distributions of the errors are quite different, especially in their tails.
Since heavier tails can be connected with higher levels of variation and  may affect the properties of  the moments,  it would be interesting to  further investigate the tails of the error distributions of high-frequency financial data.
Different tail behaviour may also be associated with different trading or marketing conditions on different days. For example,  the behaviour of the error distributions may differ on days when the whole market  or certain industrial segments such as IT are roaring.    Hence further investigations on  empirical features of this kind connecting the error distributions with practical marketing conditions can be informative for better understanding the microstructure noises in high-frequency financial data.

\section{Discussion}
\label{s5}

As the first attempt of a deconvolution approach with a contaminated slow-varying continuous process, our error density estimator method and  technical analysis in this paper assume that the distribution of the errors in the high-frequency data is symmetric. They have shown potential for revealing various  features of high-frequency financial data respecting the error distributions.   Relaxing the symmetry assumption is possible and interesting, but challenging.  The corresponding frequency domain method is under investigation.
Our analysis focuses on the univariate setting for an individual stock or trading instrument.   An interesting open question would be how to extend the frequency domain analysis to a multivariate context, which would help reveal the dependence structure between errors contaminating different stock prices.  However, the problem is significantly more difficult because of issues such as high data dimensionality and asynchronous trading records in multivariate high-frequency data.  We hope to address the problem in the future.

\section*{Acknowledgement}
We are grateful to the Editor, an Associate Editor and two referees for
their helpful suggestions. Chang was supported in part by the
Fundamental Research Funds for the Central Universities of China,
the National Natural Science Foundation of China, and
the Center of Statistical Research and the Joint Lab of Data Science and Business Intelligence at
Southwestern University of Finance and Economics. Delaigle was supported by a Future Fellowship and a Discovery Project from the Australian Research Council, and Hall was supported by a Laureate Fellowship and a Discovery Project from the Australian Research Council.
Tang acknowledges support from National Science Foundation grants.

\section*{Supplementary material}
Supplementary material available at Biometrika online includes a description of the simulation extrapolation bandwidth of \cite{DelaigleHall_2008}, additional simulation results, a construction of confidence regions for the moment estimators of Section~\ref{se:mome} and all the  proofs.

\appendix

\appendixone
\section*{Appendix}

\begin{algo}
Overview of our bandwidth selection algorithm. \label{alg:band}
\begin{tabbing}
\qquad   \quad (a) Find $f_{1}$ and $f_{_2}$ so that the relationship between $f_{1}$ and $f_{2}$ mimics that between $f_U$ and $f_{1}$,\\
\qquad   \quad\phantom{(a)}  and so that $f_{1}$ and $f_{2}$ can be accurately estimated using our data. \\
\qquad   \quad (b) In our estimator at \eqref{eq:density}, $f_U$ is estimated from $\Delta Y_{j,\ell}=Y_{t_{j+\ell}}-Y_{t_{j}} \approx U_{t_{j+\ell}}-U_{t_{j}}$, with\\
\qquad   \quad\phantom{(b)} $|t_{j+\ell}-t_j|$ being small,  $U_{t_{j+\ell}}\sim f_U$ and   $U_{t_{j}}\sim f_U$ being independent. To imitate this, \\
\qquad   \quad\phantom{(b)} for $k=1,2$ we construct versions of  $\Delta Y_{j,\ell}$, say $\Delta Y_{j,\ell}^{_k}$, such that, if $|t_{j+\ell}^{_k}-t_j^{_k}|$ is small,\\
\qquad    \quad\phantom{(b)} $\Delta Y_{j,\ell}^{_k} \approx U_{t_{j+\ell}^{_k}}^{_k}-U_{t_{j}^{_k}}^{_k}$, where $U_{t_{j+\ell}^{_k}}^{_k}\sim f_{k}$ and $U_{t_{j}^{_k}}^{_k}\sim f_{k}$ are independent.\\
\qquad   \quad (c) For $k=1, 2$, we can compute   a bandwidth $h_{k}$ and a $\xi_{k}$ for estimating $f_{k}$  using our procedure \\
\qquad \quad \phantom{(c)}  at \eqref{eq:density} applied to the data $\Delta Y_{j,\ell}^{_k}$; see \eqref{hathxi}  and \eqref{hath2}.   \\
\qquad   \quad (d) (a) suggests that $h/h_{1} \approx h_{1}/h_{2}$, so that we can take $\hat h = h_{1}^2/ h_{2}$.
\end{tabbing}
\end{algo}

\begin{algo}
Constructing $\Delta Y_{j,\ell}^{_1}$  and $t_{j}^{_1}$. \label{alg:DeltaY}
\begin{tabbing}
\qquad\quad (a) For $j=0,\ldots, n-1$, let $t_j^{_1}=(t_{j}+t_{j+1})/2$ and $Y_{1,t_j^{_1}}^{_1}=(Y_{t_{j}}+Y_{t_{j+1}})/ \surd 2$. \\
\qquad\quad\phantom{1.}~~~For $\ell> 1$ and $j=0,\ldots,n-\ell-1$, take $\Delta Y_{j,\ell}^{_1}=Y_{1,t_{j+\ell}^{_1}}^{_1}-Y_{1,t_j^{_1}}^{_1}\approx U_{1,t_{j+\ell}^{_1}}^{_1}-U_{1,t_{j}^{_1}}^{_1}$, \\
\qquad\quad\phantom{1.}~~~where $U^{_1}_{1,t_j^{_1}}= (U_{t_{j}}+U_{t_{j+1}})/ \surd 2\sim f_{1}$. Does not work for $\ell=1$ because $Y_{1,t_{j+1}^{_1}}^{_1}-Y_{1,t_j^{_1}}^{_1}$\\
\qquad\quad\phantom{1.}~~~$\approx (U_{t_{j+2}}-U_{t_{j}})/ \surd 2\not\sim f_{1}*f_{1}$. Suggests taking  $\Delta Y_{j,1}^{_1}$ as in (b).\\

\qquad\quad (b)   For $j=0,\ldots,n-2$, let $t_j^{_1}=(t_{j}+t_{j+2})/2$ and $Y_{2,t_j^{_1}}^{_1}=(Y_{t_{j}}+Y_{t_{j+2}})/\surd 2$. Take\\
\qquad\quad\phantom{1.}~~~$\Delta Y_{j,1}^{_1}=Y_{2,t_{j+1}^{_1}}^{_1}-Y_{2,t_j^{_1}}^{_1}\approx U^{_1}_{2,t_{j+1}^{_1}}-U^{_1}_{2,t_j^{_1}}$ with $U^{_1}_{2,t_j^{_1}}=(U_{t_{j}}+U_{t_{j+2}})/\surd 2\sim f_{1}$.
\end{tabbing}
\end{algo}

\begin{algo}
Constructing $\Delta Y_{j,\ell}^{_2}$ and $t_{j}^{_2}$. \label{alg:DeltaY2}
\begin{tabbing}

\qquad Step 1. For $j=0,\ldots, n-3$, let
$Y_{4,t_{j}^{_2}}^{_2}=\sum_{k=0}^{3}Y_{t_{j+k}}/2$ and $t_{j}^{_2}=\sum_{k=0}^{3}t_{j+k}/4$. \\
\qquad\phantom{1.}~~~For $\ell\geq 4$ and $j=0,\ldots, n-\ell-3$, take $\Delta Y_{j,\ell}^{_2}\equiv Y_{4,t_{j+\ell}^{_2}}^{_2}-Y_{4,t_{j}^{_2}}^{_2}\approx U^{_2}_{4,t_{j+\ell}^{_2}}-U^{_2}_{4,t_j^{_2}}$,\\
\qquad\phantom{1.}~~~where $U^{_2}_{4,t_j^{_2}}= \sum_{k=0}^{3}U_{t_{j+k}}/  2 \sim f_{2}$.  Does not work for $\ell=1,2,3$ for similar reasons\\
\qquad\phantom{1.}~~~as in Algorithm~\ref{alg:DeltaY}. For $\ell=1,2,3$, suggests taking $\Delta Y_{j,\ell}^{_2}$ as Steps 2 and 3.\\

\qquad Step 2.   For $j=0\ldots,n-6$, let $Y_{3,t_{j}^{_2}}^{_2}=\sum_{k=0,1,2,6}Y_{t_{j+k}}/2$ and $t_{j}^{_2}=\sum_{k=0,1,2,6}t_{j+k}/4$.\\
\qquad\phantom{1.}~~~For $j=0,\ldots, n-5$,  let  $Y_{2,t_{j}^{_2}}^{_2}=\sum_{k=0,1,4,5}Y_{t_{j+k}}/2$  and  $t_{j}^{_2}=\sum_{k=0,1,4,5}t_{j+k}/4$.\\
\qquad\phantom{1.}~~~For $j=0,\ldots, n-6$,  let $Y_{1,t_{j}^{_2}}^{_2}=\sum_{k=0}^{3}Y_{t_{j+2k}}/2$ and  $t_{j}^{_2}=\sum_{k=0}^{3}t_{j+2k}/4$. \\

\qquad Step 3. For $j=0,\ldots,n-9$, take
$\Delta Y_{j,3}^{_2}=Y_{3,t_{j+3}^{_2}}^{_2}-Y_{3,t_{j}^{_2}}^{_2}\approx U_{3,t_{j+3}^{_2}}^{_2}-U_{3,t_{j}^{_2}}^{_2}$, where\\
\qquad\phantom{1.}~~~$U^{_2}_{3,t_j^{_2}}= \sum_{k=0,1,2,6} U_{t_{j+k}}/  2 \sim f_{2}$.
For $j=0,\ldots,n-7$ and $\ell=1,2$, take $\Delta Y_{j,\ell}^{_2}=$\\
\qquad\phantom{1.}~~~$Y_{\ell,t_{j+\ell}^{_2}}^{_2}-Y_{\ell,t_{j}^{_2}}^{_2}\approx U_{\ell,t_{j+\ell}^{_2}}^{_2}-U_{\ell,t_{j}^{_2}}^{_2}$, where
$U^{_2}_{1,t_j^{_2}}= \sum_{k=0}^{3}U_{t_{j+2k}}/  2 \sim f_{2}$ and \\
\qquad\phantom{1.}~~~$U^{_2}_{2,t_j^{_2}}=\sum_{k=0,1,4,5}^{3}U_{t_{j+k}}/  2 \sim f_{2}$.
\end{tabbing}
\end{algo}

\end{document}